\documentclass[letterpaper, 10 pt, conference]{ieeeconf}             %
\IEEEoverridecommandlockouts        %
\usepackage{graphicx} %
\usepackage[cmex10]{amsmath} %
\usepackage{amssymb}  %
\usepackage{amsthm,thmtools,bm,dsfont}
\usepackage{hyperref}
\let\labelindent\relax
\usepackage{enumitem} \setlist[itemize,1]{leftmargin=\dimexpr 22pt}
\usepackage{cite,microtype,url}
\usepackage[bb=boondox,bbscaled=.95]{mathalfa}
\usepackage{pgfplots} 
\usepgfplotslibrary{fillbetween}
\usepackage{tikz}  
\usetikzlibrary{arrows, patterns,calc,decorations.pathmorphing,decorations.markings}

\declaretheorem[style=remark]{theorem}
\declaretheorem[style=remark]{lemma}

\declaretheorem[style=remark,qed=$\blacktriangle$]{example}
\declaretheorem[style=remark, numbered=no, ,qed=$\diamond$]{standing assumption}
\declaretheorem[style=remark]{definition}
 
\newcommand {\nn}{\nonumber}
\newcommand{\beq}{\begin{equation}}
\newcommand{\eeq}{\end{equation}}
\newcommand {\bseq}{\begin{subequations}}
\newcommand {\eseq}{\end{subequations}}
\newcommand {\bma}{\left[}
\newcommand {\ema}{\right]}

\newcommand {\Zge}{\mathbb{Z}_{+}} 		%
\newcommand {\R}{\mathbb{R}} 			%
\newcommand {\Rge}{\mathbb{R}_{+}} 	%
\newcommand {\Co}{\mathbb{C}} 			%
\renewcommand {\H}{\mathcal{H}} 			%
\renewcommand{\L}{\mathcal{L}}
\newcommand{\Li}{\mathbb{L}}
\renewcommand{\S}{\mathbb{S}}
\newcommand{\RL}{\R_{\Lambda}}
\newcommand{\SL}{\S_{\Lambda}}

\renewcommand{\Re}{\operatorname{Re}} %
\newcommand{\transpose}{\mathsf{T}} %
\newcommand{\norm}[1]{\left\lVert#1\right\rVert}

\newcommand{\spectrum}[1]{{\sigma({#1})}}
\newcommand{\overbar}[1]{\mkern 1.5mu\overline{\mkern-1.5mu#1\mkern-1.5mu}\mkern 1.5mu}
\newcommand{\one}{\delta_{-1}}
\DeclareMathOperator*{\esssup}{ess\,sup}

\newcommand{\TAC}{\textit{{IEEE} Trans. Autom. Control}}

\newcommand{\MCSS}{\textit{Math. Control, Sign. Syst.}}
\newcommand{\NOLCOS}[1]{\textit{{#1} IFAC Symp. Nonlinear Control Syst.}}

\title{\LARGE The  $\H_{\infty,p}$  norm as the differential  $\L_{2,p}$  gain of a  $p$-dominant  system}
 
\author{Alberto Padoan, Fulvio Forni, Rodolphe Sepulchre%
  \thanks{A. Padoan, F. Forni, R. Sepulchre are with Department of Engineering, University of Cambridge, Cambridge, CB2 1PZ,  UK (e-mail: \{ {\tt \small a.padoan |  f.forni | r.sepulchre \} @eng.cam.ac.uk}). The research leading to these results has received funding from the European Research Council under the Advanced ERC Grant Agreement Switchlet n. 670645  and from the Royal Society Research Grant RGS$\backslash$R1$\backslash$191308.}}

\begin{document}

\maketitle
\thispagestyle{empty}
\pagestyle{empty}

\begin{abstract} \noindent  
The differential  $\L_{2,p}$  gain of a linear, time-invariant,  $p$-dominant   system is shown to coincide with the  $\H_{\infty,p}$ norm  of its transfer function $G$,  defined as the essential supremum of the absolute value of $G$ over a vertical strip in the complex plane such that $p$ poles of $G$ lie to right of the strip.   The close analogy between the $\H_{\infty,p}$ norm and the classical $\H_{\infty}$ norm suggests that robust dominance of linear systems can be studied along the same lines as robust stability.  This property can be exploited in the analysis and design of nonlinear uncertain systems that can be decomposed as the feedback interconnection of a linear, time-invariant system with bounded gain uncertainties or nonlinearities.
\end{abstract}

\section{Introduction}

The recent paper~\cite{forni2018differential} proposes $p$-dissipativity as a generalization of the classical notion of dissipativity, with the aim of developing an interconnection theory for open $p$-dominant systems. The property of $p$-dominance formalizes the  idea   that the asymptotic behavior of a system is $p$-dimensional. The   significance   of this property for nonlinear systems analysis is apparent for small values of $p$, as the possible attractors are severely constrained  in low dimensional  systems. A $p$-dominant   system has a unique  equilibrium  point if ${p=0}$, one or several  equilibrium  points if ${p=1}$, and the simple attractors of Poincar\'e-Bendixson theorem if ${p=2}$.    In this context, $p$-dissipativity theory reformulates classical interconnection theorems of linear quadratic dissipativity theory, thus  inheriting its \textit{modus operandi} and its computational tools. The key point is that the quadratic form that characterizes the Lyapunov function or  storage is no longer required to be positive definite. Instead, it is required to have a fixed inertia, with $p$ negative eigenvalues and ${n-p}$ positive eigenvalues, where $n$ is the dimension of the system.

 A   notion of $\L_{2,p}$-gain can be defined for a $p$-dominant system with rate $\lambda$ using the differential  dissipation inequality
\beq \label{eq:dissipativity-inequality-sg}
\!\!
\begingroup %
\setlength\arraycolsep{1.5pt}
\bma
\begin{array}{c}
\delta \dot{x} \\
\delta x
\end{array}
\ema^{\transpose}
\!\!
\bma
\begin{array}{cc}
0 & P \\
P &  2\lambda P + \varepsilon I
\end{array}
\ema
\!
\bma
\begin{array}{c}
\delta \dot{x} \\
\delta x
\end{array}
\ema 
\! \le \!
\bma
\begin{array}{c}
\delta y \\
\delta  u
\end{array}
\ema^{\transpose}
\!\!
\bma
\begin{array}{cc}
-I & 0 \\
0 &  \gamma^2 I
\end{array}
\ema
\!
\bma
\begin{array}{c}
\delta y \\
\delta u
\end{array}
\ema  \! , \!
\endgroup
\eeq
with ${P\in\R^{n\times n}}$ a symmetric matrix with $p$ negative eigenvalues and ${n-p}$ positive eigenvalues. 
For $p=0$, the differential dissipation inequality~\eqref{eq:dissipativity-inequality-sg} simply means that the classical ${\L_2}$ gain of the system does not exceed $\gamma$.  By the KYP lemma~\cite{yakubovich1962solution,kalman1963lyapunov,brogliato2007dissipative}, $\gamma$ also coincides with the classical $\H_{\infty}$ norm of the transfer function of the system. 
Similarly,  the  $\L_{2,p}$-gain  of a finite-dimensional, linear, time-invariant, $p$-dominant system with rate $\lambda$ with transfer function $G$ can be expressed as 
\beq  \label{eq:intro-H-inf-lambda} 
 \esssup_{\omega \in \R} |G(i\omega -\lambda)|, 
\eeq
where ${G(s-\lambda)}$ has $p$ poles in the open right half-plane and ${n-p}$ poles in the open left half-plane.
This raises the question of computing the  $\L_{2,p}$-gain  of a system through~\eqref{eq:intro-H-inf-lambda} as a generalization of the classical $\H_\infty$ norm.

The goal of the paper is to outline an  $\H_{\infty,p}$   theory geared towards $p$-dominance that closely parallels classical $\H_\infty$ theory.  The $\H_{\infty,p}$ norm for functions defined on a vertical strip  is shown to be the system norm induced by   the unique bounded operator defined by a transfer function with $p$ poles to the right of its region of convergence and $n-p$ poles to the left of its region of convergence.  The paper   emphasizes that most usual properties of the classical  $\H_\infty$ norm carry over to the $\H_{\infty,p}$ norm.  The motivation   is to use the $\H_{\infty,p}$ norm for robustness and performance analysis of $p$-dominant systems in the same way as one uses the $\H_{\infty}$ norm for stable systems.

The paper is organized as follows. 
Section~\ref{sec:main-results} introduces  Hardy spaces   on a vertical strip. Section~\ref{sec:operator}  shows that the $\H_\infty$ norm for  Hardy spaces   on a vertical strip can be interpreted as a norm induced by a multiplication operator and by a convolution operator on the whole real line.  Section~\ref{sec:dominance}  illustrates some connections between  Hardy spaces  on a vertical strip and $p$-dominance theory~\cite{forni2018differential,felix2018analysis,padoan2019feedback}. 
Section~\ref{sec:example} provides an illustrative example of robust $p$-dominance analysis. 
Section~\ref{sec:conclusion} concludes the paper with some final remarks and future research directions.  The appendix provides additional background material on the bilateral Laplace transform. 
The proofs are omitted for reasons of space.  

\textit{Notation}: 
$\R$ and $\Co$ denote the set of real numbers and the set of complex numbers, respectively. 
$\Zge$ and $\Rge$ denote the set of non-negative integer numbers and the set of non-negative real numbers, respectively.  $i$ denotes the imaginary unit and $i\R$ denotes the set of complex numbers with zero real part. $\partial{S}$ denotes the boundary of the set $S$. 
$I$ denotes the identity matrix. 
$\spectrum{A}$ denotes the spectrum of the matrix ${A\in\Co^{n \times n}}$.  $M^{\transpose}$ and $M^{*}$ denote the transpose and the conjugate transpose of the matrix ${M\in\Co^{l \times m}}$, respectively. 
${|\cdot|}$ denotes the standard Euclidean norm on~${\Co^n}$.

\section{Hardy  spaces   on a vertical strip} \label{sec:main-results}

This section introduces  Hardy spaces   of functions on a vertical strip. 
Let  ${q\in\Zge}$, with ${q\ge 1}$,    let~${\Lambda = (\underline{\lambda}, \overbar{\lambda})}$ be an open interval, with ${ -\infty \le \underline{\lambda} < \overbar{\lambda} \le \infty}$, and let
$\SL  = \left\{ \, s\in \Co \, :\,   \Re(s) \in  -\Lambda  \,\right\} $, 
with ${-\Lambda = \left\{ \, - \lambda \, :\,   \lambda \in \Lambda \,\right\} }$.

\begin{definition} \label{def:Hardy-strip}
${\H_q(\SL)}$ is the set of all analytic functions\footnote{Lebesgue integration is used throughout this work. Functions that are equal except for a set of measure zero are identified. Conditions imposed on a function are understood in the sense of being valid for all points of the domain of the function except for a set of measure zero.} ${f: \SL \to\Co^n}$ such that 
\beq 
\esssup_{\lambda \in \Lambda}{\left(\int_{-\infty}^{\infty} |f(  -\lambda  +i \omega)|^{q} \frac{d\omega}{2 \pi} \right) < \infty}.
\eeq
${\H_\infty(\SL)}$ is the set of all analytic functions ${f: \SL \to\Co^n}$ such that
\beq
{\displaystyle \esssup_{s\in \SL } |f(s)| < \infty}.
\eeq
The \textit{${\H_q(\SL)}$ norm} of a function ${f\in \H_q(\SL)}$ is defined as
\beq \label{eq:Hq-norm-lambda-strip}
\scalebox{.87}{$
\norm{f}_{\H_q(\SL)} \!  = \! 
 \left\{ \! 
\begingroup %
\setlength\arraycolsep{1.2pt}
   	\begin{array}{ll}
     \displaystyle \esssup_{\lambda \in \Lambda}{\left(\int_{-\infty}^{\infty} |f( -\lambda   +i \omega)|^{q} \frac{d\omega}{2 \pi} \right)^{\frac{1}{q}}}  \!  , & \text{for } 1 \le q < \infty  , \\
     \displaystyle \esssup_{s\in  \SL } |f(s)| , & \text{for } q    = \infty  .
     \end{array}
  \endgroup %
    \right.
    $}
\eeq 
\end{definition}

\noindent  
${\H_q(\SL)}$   is a linear space,  with scalar product and sum defined in the standard fashion.   It is therefore  referred to as  a   \textit{Hardy space on the vertical strip $\SL$}, as it possesses many of the nice properties of classical Hardy spaces.  

\begin{theorem} \label{thm:Banach}
${\H_q(\SL)}$ is a Banach space for ${1\le q\le \infty}$.
\end{theorem}

\begin{theorem} \label{thm:hardy-strip-max-first} 
Let ${f\in \H_\infty(\SL)}$.  Assume $f$ is continuous and bounded  on $\partial \SL$.   Then
\beq 
\norm{f}_{\H_{\infty}(\SL)}  
= \displaystyle \sup_{s \in \partial \SL} |f(s)| .
\eeq
\end{theorem}

\noindent
Theorem~\ref{thm:hardy-strip-max-first} establishes a maximum modulus theorem for functions in ${\H_\infty(\SL)}$:  the norm of a function ${f\in \H_\infty(\SL)}$ can be computed by only considering the behavior of $f$ on the \textit{boundary} of the strip ${\SL}$, provided that $f$ is continuous and bounded therein. Thus the following standing assumption is made in order to apply Theorem~\ref{thm:hardy-strip-max-first} throughout the paper.

\begin{standing assumption}  
Every function ${f\in \H_\infty(\SL)}$ is continuous  and bounded on $\partial \SL$.  
\end{standing assumption}

The classical $\H_\infty$ norm of a function is tightly connected to the $\L_\infty$ norm of the corresponding boundary function defined on the imaginary axis~\cite[p.7]{partington2004linear}. We now show that a similar property holds for the $\H_\infty(\SL)$ norm. Given ${\lambda \in \R}$ 
let %
${ \Li_{\lambda} = \left\{ \, s\in \Co \, :\,  \Re(s) =  -\lambda   \,\right\} }.$

\begin{definition}
${\L_q(\Li_{\lambda})}$ is   the set of all measurable functions ${f: \Co \to\Co^n}$ such that
\beq 
\int_{-\infty}^{\infty} |f(  -\lambda   + i\omega)|^{q} \frac{d\omega}{2 \pi}  < \infty.
\eeq
${\L_\infty(\Li_{\lambda})}$ is   the set of all measurable functions ${f: \Co \to\Co^n}$ such that
\beq
{\displaystyle \esssup_{\omega\in \R} |f( -\lambda   + i\omega)| < \infty}.
\eeq
The \textit{${ \L_q(\Li_{\lambda})}$ norm} of a function ${f\in \L_q(\Li_{\lambda})}$ is defined as
\beq  \label{eq:Lq-norm-frequency-lambda}
\norm{f}_{ \L_q(\Li_{\lambda})} \!  = \! 
 \left\{ \! 
\begingroup %
\setlength\arraycolsep{1.4pt}
   	\begin{array}{ll}
     \displaystyle {\left( \int_{-\infty}^{\infty} |f( -\lambda   +i\omega)|^{q} \frac{d\omega}{2 \pi} \right)^{\frac{1}{q}}} \!  , & \text{for } 1 \le q < \infty  , \\
     \displaystyle \esssup_{\omega\in \R} |f( -\lambda  +i\omega)| , & \text{for } q    = \infty ,
     \end{array}
  \endgroup %
    \right.
\eeq 
\end{definition}

\noindent
The norm~\eqref{eq:Lq-norm-frequency-lambda} induces a Banach space structure on the set ${\L_q(\Li_{\lambda})}$ for ${1 \le q \le \infty}$~\cite[p.19]{devore1993constructive}. For ${q=2}$, it coincides with the norm induced by the inner product
\beq  \label{eq:L2-inner-product-frequency-lambda}
\langle f, g\rangle_{\L_2(\Li_{\lambda})} = \frac{1}{2\pi}\int_{-\infty}^{\infty} f(-\lambda+\omega)^{*} g(-\lambda+\omega) d\omega .
\eeq
${\L_2(\Li_{\lambda})}$ is therefore as a Hilbert space, which admits the (orthogonal direct sum) decomposition 
\beq \label{eq:decomposition-Lilambda}
{\L_2(\Li_{\lambda})= \H_2(\S_{\Lambda_{-}}) \oplus \H_2(\S_{\Lambda_{+}}) ,} 
\eeq
with ${\Lambda_{-} =  (\lambda,\infty)  }$ and ${\Lambda_{+} =  (-\infty, \lambda)  }$, in which the orthogonality condition ${\langle f, g\rangle_{\L_2(\Li_{\lambda}) }= 0}$ holds for every ${f\in \H_2(\S_{\Lambda_{-}})}$ and ${g\in  \H_2(\S_{\Lambda_{+}}) }$.

We are now ready to connect ${\H_{\infty}(\SL)} $ and $\L_\infty(\Li_{\lambda})$ norms.

\begin{theorem} \label{thm:hardy-strip-max}
Under the assumption of Theorem~\ref{thm:hardy-strip-max-first},
\beq 
\norm{f}_{\H_{\infty}(\SL)}  
= \max\{ 
	\norm{f}_{\L_{\infty}(\Li_{\underline{\lambda}})}, 
	\norm{f}_{\L_{\infty}(\Li_{\overbar{\lambda}})}
	\} .
\eeq
\end{theorem}

\noindent
Theorem~\ref{thm:hardy-strip-max} is consistent with the classical ``limit'' cases.   For  ${\Lambda = \Rge}$ (${-\Lambda = \Rge}$)  the strip $\SL$ is the  open half-plane to the right (left) of the imaginary axis $i\R$  and the ${\H_\infty(\SL)}$ norm reduces to the norm
\beq \label{eq:Hq-norm-lambda1}
\norm{f}_{\H_\infty(\SL)} = \displaystyle \esssup_{s\in \SL } |f(s)| = \norm{f}_{\L_{\infty}(i\R)} .
\eeq  
For ${\Lambda =(\underline{\lambda}, \overbar{\lambda})}$, with ${\underline{\lambda} \to \lambda^-}$  and ${ \overbar{\lambda} \to \lambda^+}$,  the strip $\SL$  tends to the vertical line $\Li_{\lambda}$ and the ${\H_\infty(\SL)}$ norm reduces to the norm
\beq \label{eq:Hq-norm-lambda3}
\norm{f}_{\H_\infty(\SL)} = \displaystyle \esssup_{s\in \Li_{\lambda} } |f(s)| = \norm{f}_{\L_{\infty}(\Li_{\lambda})}.
\eeq

\section{The ${\H_\infty(\SL)}$ norm as an induced norm} \label{sec:operator}

A classical result of $\H_\infty$ theory is that the norm induced by the multiplication operator associated with a function ${G \in \L_{\infty}(i\R)}$ coincides with the $\L_{\infty}$ norm of $G$~\cite[p.100]{zhou1996robust}. We now show that a similar result holds if ${G\in\H_{\infty}(\SL)}$. 

\begin{definition}
The \textit{multiplication operator} associated with the function ${G\in \H_{\infty}(\SL)}$ is defined as
\beq  \label{eq:multiplication-operator}
M_G: \H_2(\SL) \to \H_2(\SL), \ U \mapsto G U,
\eeq
and the corresponding $\H_{2}(\SL)$ induced norm is defined as 
\beq \label{eq:induced-norm-multiplication-operator}
\norm{M_G}_{\H_{2}(\SL)} 
=    \sup_{\substack{U \in \H_{2}(\SL) \\ \norm{U}_{\H_{2}(\SL)} = 1}}   \norm{GU}_{\H_{2}(\SL)} .
\eeq
\end{definition}

\begin{theorem} \label{thm:multiplication-operator-norm} 
Let ${G\in \H_{\infty}(\SL)}$ and consider the    multiplication operator~\eqref{eq:multiplication-operator}. Then  ${\norm{M_G}_{\H_{2}(\SL)} = \norm{G}_{\H_\infty(\SL) }.}$
\end{theorem}

\noindent
The ${\H_{\infty}(\SL)} $ norm can be also characterized as the norm induced by the convolution operator associated with a continuous-time, single-input, single-output,  linear, time-invariant system described by the equations 
\beq \label{eq:system-linear-MIMO}
\quad \dot{x} = Ax+Bu, \quad y=Cx+Du,
\eeq
with ${x(t)\in\R^n}$, ${u(t)\in\R}$, ${y(t)\in\R}$, ${A\in\R^{n \times n}}$, ${B\in\R^{n}}$, ${C\in\R^{1\times n}}$ and ${D\in\R}$ constant matrices, and transfer function ${G(s) = C(sI-A)^{-1} B +D}$.

If system~\eqref{eq:system-linear-MIMO} has no eigenvalues in $\SL$, then a (unique) bounded convolution operator can be associated with the system by defining its impulse response as 
\beq \label{eq:unstable-impulse-response}
g(t) 
=
\begin{cases}
	C_{+} e^{A_+ t} B_+ , & \text{ for } t >0 , \\
	C_{-} e^{A_- t} B_- , & \text{ for } t \le 0 ,
\end{cases}
\eeq
in which, upon a possible coordinates change,
\beq 
A =
\bma
\begin{array}{cc}
A_+ & 0 \\
0 & A_-
\end{array}
\ema, \
B =
\bma
\begin{array}{c}
B_+ \\
B_-
\end{array}
\ema, \
C^{\transpose} =
\bma
\begin{array}{c}
C_+^{\transpose} \\
C_-^{\transpose}
\end{array}
\ema, 
\eeq
with  ${\spectrum{A_+} \subset \S_{\Lambda_{+}}}$ and ${\spectrum{A_-} \subset \S_{\Lambda_{-}}}$ for  
${{\Lambda_{+}} = (-\infty,\underline{\lambda}) }$ and ${{\Lambda_{-}} = (\overbar{\lambda},\infty) }$,  respectively. The impulse response~\eqref{eq:unstable-impulse-response} is uniquely defined by the inverse  (bilateral) Laplace transform of $G$ in its region of convergence $\SL$.   Conversely,  the transfer function of system~\eqref{eq:system-linear-MIMO} coincides with the Laplace transform of its impulse response, \textit{i.e.} ${ G(s) = \L\{g\}(s)}$. 

\begin{definition}
\textit{${\L_{q} (\RL)}$}  is defined as the set of all measurable functions ${f: \R \to\Co^n}$ such that
\beq 
\esssup_{\lambda \in \Lambda} \left( \displaystyle  \int_{-\infty}^{\infty}   e^{q\lambda t}   | f(t) |^{q}  dt  \right) <\infty .
\eeq
The  \textit{${\L_{q}(\RL)}$ norm} of a function ${f\in \L_{q}(\RL)}$ is defined as
\beq  \label{eq:L2-uniformly-weighted-norm}
\norm{f}_{\L_{q}(\RL)} = \esssup_{\lambda \in \Lambda} \left( \displaystyle  \int_{-\infty}^{\infty}  
 e^{q\lambda t}  | f(t) |^{q}  dt   \right)^{\frac{1}{q}}  .
\eeq
\end{definition} 

\noindent
\begin{definition}
The \textit{convolution operator} associated with system~\eqref{eq:system-linear-MIMO} is defined as\footnote{The  same symbol is used for the convolution operator associated with a system and the corresponding transfer function. Context determines which is meant.}
\beq \label{eq:unstable-convolution-operator}
\!\!\scalebox{.89}{$ 
G:  \L_2(\RL) \to \L_2(\RL) ,  \, u \mapsto \displaystyle  \int_{-\infty}^{\infty} g(t-\tau) u(\tau)d\tau + Du(t), 
$} \!
\eeq
and the corresponding induced $\L_2(\RL)$ norm is defined as 
\beq \label{eq:unstable-induced-norm}
\norm{G}_{\L_2(\RL) } 
=    \sup_{\substack{u \in \L_2(\RL) \\ \norm{u}_{\L_2(\RL)} = 1}}   \norm{G u}_{\L_2(\RL)} .
\eeq
\end{definition}

\begin{theorem} \label{thm:convolution-operator-norm}
Consider system~\eqref{eq:system-linear-MIMO} and the  associated   convolution operator~\eqref{eq:unstable-convolution-operator}. Then ${\norm{G}_{\L_2(\RL)}  = \norm{G}_{\H_\infty(\SL) }}.$
\end{theorem}

 \noindent 
Theorem~\ref{thm:convolution-operator-norm} establishes that the $\L_2(\RL)$  norm induced by the   convolution operator associated with system~\eqref{eq:system-linear-MIMO}   coincides with the $\H_{\infty}(\SL)$ norm of the transfer function  of system~\eqref{eq:system-linear-MIMO}.    Thus the $\H_{\infty}(\SL)$ norm of a transfer function can be interpreted as the \textit{gain} of the corresponding system, in analogy with classical $\H_{\infty}$ theory~\cite{zhou1996robust}.

\section{The $\H_{\infty}(\SL)$ space and dominant systems} \label{sec:dominance}

\subsection{The differential  $\L_{2,p}$  gain of a  $p$-dominant  system}

The discussion above is of interest because of its applications to $p$-dominance theory~\cite{forni2018differential}. In what follows we summarize relevant definitions and properties.
Consider a continuous-time, nonlinear, time-invariant system and its linearization described by the equations
\bseq \label{eq:system-nonlinear-prolonged}
\begin{align}
\dot{x} &= f(x) +Bu, \quad \quad \quad \ \   y = \ Cx+ D  u,    \label{eq:system-nonlinear} \\
\delta\dot{x} &= \partial f(x) \delta x + B \delta u,   \quad \delta y = C \delta y + D \delta u,   \label{eq:system-nonlinear-linearization}
\end{align}
\eseq
in which ${x(t)\in\R^n}$, ${u(t)\in\R^m}$, ${y(t)\in\R^m}$, 
${f:\R^{n} \to \R^{n}}$ is a continuously differentiable\footnote{This assumption simplifies the exposition. Analogous considerations can be performed requiring only Lipschitz continuity.} vector field, ${B\in\R^{m\times n}}$, ${C\in\R^{m\times n}}$ and ${D\in\R^{m\times m}}$ are constant matrices, ${\delta x(t)\in\R^n}$, ${\delta u(t)\in\R^m}$, ${\delta y(t)\in\R^m}$ (identified with the respective tangent spaces), and ${\partial f}$ is the Jacobian of the vector field $f$.

\begin{definition} \cite{forni2018differential} \label{def:dominant-linear} 
For ${u=0}$, system~\eqref{eq:system-nonlinear} is \textit{$p$-dominant with rate ${\lambda \in \Rge}$}  if   there exist ${\varepsilon \in \Rge}$  and a symmetric matrix ${P\in\R^{n\times n}}$, with inertia\footnote{\!\!
The inertia of the matrix ${A \in \R^{n \times n}}$ is defined as  $(\nu,\delta, \pi)$, where $\nu$ is the number of eigenvalues of $A$ in the open left half-plane, $\delta$ is the number of eigenvalues of $A$ on the imaginary axis, and $\pi$ is the number of eigenvalues of $A$ in the open right half-plane, respectively.  
}  $(p, 0, n-p)$,   such that the conic constraint
\beq \label{eq:system-conic-constraint-dominance}
\bma
\begin{array}{c}
\delta \dot{x} \\
\delta x
\end{array}
\ema^{\transpose}
\bma
\begin{array}{cc}
0 & P \\
P &  2\lambda P + \varepsilon I
\end{array}
\ema
\bma
\begin{array}{c}
\delta \dot{x} \\
\delta x
\end{array}
\ema \le 0 
\eeq 
holds along the solutions of the prolonged system~\eqref{eq:system-nonlinear-prolonged}. The property is strict if ${\varepsilon >0}$.
\end{definition}

\begin{definition} \label{def:pgain}
The system~\eqref{eq:system-nonlinear} is said to have \textit{(finite differential) ${\L_{2,p}}$-gain (from $u$ to $y$) less than ${\gamma \in \Rge}$ with rate ${\lambda \in \Rge}$} if there exist $\varepsilon \in \Rge$ and a symmetric matrix ${P\in\R^{n\times n}}$, with inertia   $(p, 0, n-p)$,   such that the conic constraint
\beq \label{eq:system-linear-MIMO-conic-IO}
\!\!
\scalebox{.925}{$
\begingroup %
\setlength\arraycolsep{1.5pt}
\bma
\begin{array}{c}
\delta \dot{x} \\
\delta  x
\end{array}
\ema^{\transpose} \!
\bma
\begin{array}{cc}
0 & P \\
P &  2\lambda P + \varepsilon I
\end{array}
\ema \!
\bma
\begin{array}{c}
\delta  \dot{x} \\
\delta  x
\end{array}
\ema 
\le
\bma
\begin{array}{c}
\delta y \\
\delta  u
\end{array}
\ema^{\transpose} \!
\bma
	\begin{array}{cc}
	-I & 0 \\
	 0 & \gamma^2 I
	\end{array}
	\ema \!
\bma
\begin{array}{c}
\delta y \\
\delta  u
\end{array}
\ema  \! 
\endgroup
$}
\eeq
holds along the solutions of the prolonged system~\eqref{eq:system-nonlinear-prolonged}. The \textit{(differential) $\L_{2,p}$-gain of system~\eqref{eq:system-nonlinear}  (from $u$ to $y$) with rate $\lambda$} is defined as ${ \gamma_{\lambda} = \inf \left\{ \gamma \in \Rge  :  \text{\eqref{eq:system-linear-MIMO-conic-IO} holds} \right\}}.$ The properties are strict if $\varepsilon >0$. 
\end{definition}

 \noindent 
The property of $p$-dominance strongly constrains the asymptotic behavior of a system,  as clarified by the next theorem.

\begin{theorem} \cite{forni2018differential} \label{thm:asymptotic}
Assume system~\eqref{eq:system-nonlinear} is strictly $p$-dominant with rate ${\lambda \in \Rge}$  and let ${u=0}$.  Then every bounded solution of~\eqref{eq:system-nonlinear} converges asymptotically~to 
\begin{itemize}[leftmargin=\dimexpr 14pt]
\item the unique equilibrium  point if ${p=0}$,
\item a (possibly non-unique) equilibrium  point if ${p=1}$,
\item an equilibrium  point, a set of equilibrium  points and their connected arcs or a limit cycle if ${p=2}$.
\end{itemize}
\end{theorem} 
 
 \noindent 
The ${\L_{2,p}}$-gain can be used to establish $p$-dominance of an interconnected system, thus extending classical small-gain conditions~\cite{desoer1975feedback}.

\begin{theorem}[Small-gain theorem for $p$-dominance] \label{thm:small-gain}
Let $\Sigma_i$ be a system with input ${u_i\in\R^{m_i}}$, output ${y_i\in\R^{m_i}}$, and (differential) $\L_{2,p_i}$-gain less than ${\gamma_i\in\Rge}$ with rate ${\lambda \in \Rge}$, with ${i \in \{1,2\}}$. Then the closed-loop system $\Sigma$ defined by the negative feedback interconnection equations  ${u_1 = -y_2}$ and ${u_2 = y_1}$  is $(p_1+p_2)$-dominant with rate $\lambda$ if ${\gamma_1 \gamma_2 <1}$. 
\end{theorem}

\subsection{The  ${\H_{\infty,p}(\SL)}$  norm as the differential  $\L_{2,p}$  gain  }

For a linear, time-invariant system~\eqref{eq:system-linear-MIMO} the conic constraint~\eqref{eq:system-linear-MIMO-conic-IO} holds  along the solutions of the system if and only if  there exist ${\varepsilon \in \Rge}$ and a symmetric matrix ${P\in\R^{n\times n}}$, with inertia   $(p, 0, n-p)$,   which solve  the linear matrix inequality
\beq \label{eq:linear-LMI}
\scalebox{.94}{$
\begingroup %
\setlength\arraycolsep{1pt}
\bma
\begin{array}{cc}
A^{\transpose} P + P A + 2\lambda P + \varepsilon I + C^{\transpose} C 
&  PB   + C^{\transpose}D \\
B^{\transpose}P  + D^{\transpose}C
&  -\gamma^2I+D^{\transpose}D
\end{array}
\ema
\endgroup
$} \le 0  .
\eeq
In particular, system~\eqref{eq:system-linear-MIMO} is $p$-dominant if and only if  there exist  $\varepsilon \in \Rge$  and a symmetric matrix ${P\in\R^{n\times n}}$, with inertia   $(p, 0, n-p)$, such that   
\beq \label{eq:linear-LMI-pdominance}
A^{\transpose} P + P A + 2\lambda P + \varepsilon I  \le 0.
\eeq 
The inertia constraint in~\eqref{eq:linear-LMI} and~\eqref{eq:linear-LMI-pdominance} entails that the transfer function $G$ has $p$ poles to the right of the line $\Li_{\lambda}$ and $n-p$ poles to the left of the line $\Li_{\lambda}$. If $G$  satisfies this property for every ${\lambda\in \Lambda}$, then  ${G\in \H_{\infty}(\S_{\Lambda})}$.   The converse is also true: if  ${G\in \H_{\infty}(\S_{\Lambda})}$   is rational, then $G$ can be realized by a $p$-dominant system with rate ${\lambda \in \Lambda}$.  The Hardy space ${\H_\infty(\SL)}$ and $p$-dominant systems are strongly related.   Thus, it is convenient to introduce the notation ${\H_{\infty,p}(\SL)}$ for the subspace of all functions in ${  \H_{\infty}(\SL)}$ with $p$ poles in the open half-plane to the right of $\SL$. 
The next result clarifies the  interplay between  the $\H_{\infty,p}(\SL)$ norm of a transfer function the \textit{$\L_{2,p}(\RL)$-gain} of the corresponding system, defined as 
\beq 
\gamma_{\Lambda} = \sup_{\lambda \in \Lambda} \gamma_{\lambda} ,
\eeq
with ${\gamma_{\lambda}  \in \Rge}$ the $\L_{2,p}$-gain of system~\eqref{eq:system-linear-MIMO} with rate ${\lambda \in \Lambda}$.

\begin{theorem} \label{thm:pgain}
Let   ${\gamma_{\Lambda}  \in \Rge}$   be the  $\L_{2,p}(\RL)$-gain   of system~\eqref{eq:system-linear-MIMO}. Then 
${\gamma_{\Lambda}  =  \norm{G}_{\H_{\infty,p}(\S_{\Lambda})}.  }$
\end{theorem}

\noindent
Theorem~\ref{thm:pgain}  connects the state-space notion of  $\L_{2,p}(\RL)$-gain   of a system to the frequency-domain notion of $\H_{\infty,p}(\SL)$ norm of the corresponding transfer function.  This   opens the way to robust $p$-dominance analysis and design through standard state-space and frequency-domain tools, including linear matrix inequalities, Riccati equations, Bode diagrams and Nyquist diagrams.

\subsection{Computation of the $\H_{\infty}(\SL)$ norm}

By Theorem~\ref{thm:hardy-strip-max}, the norm of a  function  ${G\in\H_{\infty}(\SL)}$ can be computed as the maximum between the $\L_{\infty}(\Li_{\underline{\lambda}})$ norm and the $\L_{\infty}(\Li_{\overbar{\lambda}})$ norm of $G$. These norms, in turn, coincide with $\L_{\infty}(i\R)$  norm of the $\lambda$-shifted transfer function ${G_{\lambda} : s \mapsto G(s-\lambda)}$ for  ${\lambda = \underline{\lambda}}$ and ${\lambda = \overbar{\lambda}}$,   respectively. As a result, the  $\H_{\infty}(\SL)$ norm of $G$ can be computed by first considering the $\lambda$-shifted transfer function $G_{\lambda}$ for  ${\lambda = \underline{\lambda}}$ and ${\lambda = \overbar{\lambda}}$,  then  computing their $\L_{\infty}(i\R)$ norm, and finally taking the maximum between the two values.

If $G$ is rational, these computations can be efficiently performed using established state-space methods~\cite{zhou1996robust}. For example, the  $\L_{\infty}(\Li_{\lambda})$   norm of the transfer function $G$ can be computed via a bisection algorithm based on testing if the condition
\beq \label{eq:bisection-condition}
\norm{G}_{\L_{\infty}   (\Li_{\lambda}) }   < \gamma
\eeq
holds for a given constant ${\gamma\in\Rge}$. The test can be performed by solving the linear matrix inequality~\eqref{eq:linear-LMI} based on the fact that the linear, time-invariant system~\eqref{eq:system-linear-MIMO} has $\L_{2,p}$-gain less than $\gamma$ if and only if the linear matrix inequality~\eqref{eq:linear-LMI} admits  a solution.   By Theorem~\ref{thm:pgain}, this means that if system~\eqref{eq:system-linear-MIMO} is a realization of the transfer function $G$, then the condition~\eqref{eq:bisection-condition} holds if and only if the linear matrix inequality~\eqref{eq:linear-LMI} admits a solution.  

A similar bisection algorithm can be devised by checking iteratively if the Hamiltonian matrix
\beq  \label{eq:bisection-hamiltonian}
 \!  \! 
\scalebox{0.9}{$
\begingroup %
\setlength\arraycolsep{1.5pt}
H_\gamma \! = \! 
\bma
\begin{array}{cc}
A+\lambda I - BR^{-1}D^{\transpose}C
& -\gamma BR^{-1}B^{\transpose}  \\
-\gamma C^{\transpose}S^{-1}C
&   -A^{\transpose}-\lambda I - C^{\transpose}DR^{-1}B^{\transpose}
\end{array}
\ema  
 \endgroup $},
\eeq
with ${\gamma\in\Rge}$, ${R=D^{\transpose} D - \gamma^{2}I}$ and ${S=D D^{\transpose} - \gamma^{2}I}$, has eigenvalues on the imaginary axis. This is a consequence of~\cite[Theorem 1]{boyd1989bisection}, which we recall below for completeness.

\begin{theorem}
Consider system~\eqref{eq:system-linear-MIMO}. Assume ${\spectrum{A} \cap i\R = \emptyset}$, ${\gamma\in\Rge}$ is not a singular value of $D$, and ${\omega_0 \in \R}$. Then $\gamma$ is a singular value of $G(i\omega_0)$ if and only if $(H_{\gamma} - i\omega_0 I)$ is singular.
\end{theorem}

\noindent
An estimate of the  $\L_{\infty}(\Li_{\lambda})$    norm can be also obtained as 
\beq
{\norm{G}_{\L_{\infty}   (\Li_{\lambda})   } \approx \max_{1\le k \le \nu}} |G_{\lambda}(i\omega_k )|, 
\eeq
provided that the grid of frequency points ${\{ \omega_1, \ldots, \omega_{\nu}\}}$ is sufficiently fine.  

In principle, the   $\L_{\infty}(\Li_{\lambda})$   norm of a transfer function $G$ can be also obtained graphically, as the distance in the complex plane from the origin to the farthest point on the Nyquist diagram of the $\lambda$-shifted transfer function $G_{\lambda}$ or as the peak value of the Bode diagram of the magnitude of the $\lambda$-shifted transfer function $G_{\lambda}$. Finally, the   $\L_{\infty}(\Li_{\lambda})$   norm also coincides with the essential supremum of the restriction of a transfer function along the axis   $\Li_{\lambda}$.

\section{An illustrative example} \label{sec:example}

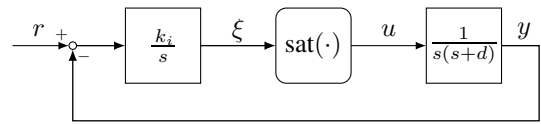
\begin{figure}[h!]
\centering 
\begin{tikzpicture}[black]
\node at (-1,3.2) {$\xi$};
\draw  (1.5,3.5) rectangle (2.5,2.5);
\node at (2,3) {$\frac{1}{s(s+d)}$};
\draw  (-2.5,3.5) rectangle (-1.5,2.5);
\node at (-2,3) {$\frac{\,k_i}{s}$};
\draw[-latex, line width = .5 pt] (-3.15,3) -- (-2.7,3) -- (-2.5,3);
\draw[rounded corners] (-0.5,3.5) rectangle (0.5,2.5);
\node at (0,3) {sat$(\cdot)$};
\node at (1,3.2) {$u$};
\node at (2.8,3.2) {$y$};
\draw [-latex, line width = .5 pt](2.5,3) -- (3,3) -- (3,2) -- (-3.2,2) -- (-3.2,2.95);
\draw [-latex, line width = .5 pt](-1.5,3) -- (-0.5,3);
\draw [-latex](0.5,3) -- (1.5,3);
\draw [-latex](-4,3) -- (-3.25,3);
\draw[fill = white] (-3.2,3) circle [radius=0.05]; 
\node at (-3.35,3.15) {\tiny $+$};
\node at (-3.05,2.85) {\tiny $-$};
\node at (-3.65,3.2) {$r$};
\end{tikzpicture}
\centering
\caption{The system~\eqref{eq:system-MSD}.}
\label{fig:feedback_interconnection_example_CDC}
\end{figure}

Consider a one-degree-of-freedom mechanical system subject to saturated integral control described by the equations 
\beq \label{eq:system-MSD}
\ddot{y} + d \dot{y}= u , \ \,
\dot{\xi} = k_i(r-y),  \ \,
u = \text{sat}\left(\xi \right),
\eeq
in which ${y(t) \in \R}$ is the position of the point mass,  ${\xi(t) \in \R}$ is the integrator variable,  ${u(t)\in \R}$ is the control input, ${r(t)\in \R}$ is the reference signal, ${d \in \Rge}$ is the damping coefficient, ${k_i \in \R}$ is the integral gain, and ${\text{sat}:\R\to\R}$ is defined as ${\text{sat}(y) = \min(\max(y, -1), 1)}$ for every ${y\in\R}$,  respectively.   The setup is illustrated in Fig.~\ref{fig:feedback_interconnection_example_CDC}.

The dominance properties of the closed-loop system~\eqref{eq:system-MSD} can be modulated through the integral gain $k_i$. By the circle criterion for $p$-dominance~\cite{felix2018analysis}, for ${r = 0}$ and for ${k_i}$ sufficiently small the system is strictly $2$-dominant with rate ${\lambda \in \Lambda}$ for every $\Lambda \subset (0,d)$. By Theorem~\ref{thm:asymptotic}, the behavior of the closed-loop system is therefore oscillatory, since its solutions are bounded and the unique equilibrium at the origin is unstable.

These conclusions have been drawn by neglecting actuator dynamics, which can be modeled in first approximation as a first order lag with transfer function 
\beq 
H(s) = \frac{1}{1+s\tau}, \qquad \tau \in \Rge.
\eeq
Actuator dynamics are indeed negligible provided they are sufficiently fast. This is well-known in the case of stability. The theory developed in the present paper allows one to extend this principle to switching and oscillatory regimes.

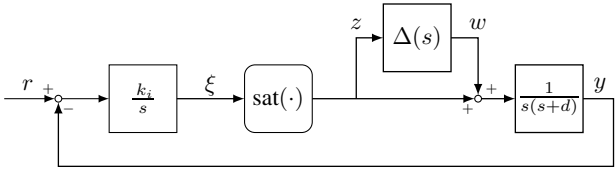
\begin{figure}[t!]
\centering 
\begin{tikzpicture}[black, scale=0.9, every node/.style={transform shape}]
 \node at (3,4.1) {$w$};
 \node at (1.2,4.1) {$z$};
\draw[fill = white] (3,3) circle [radius=0.05]; 
\draw [-latex, line width = .5 pt] (1.6,4.4) rectangle (2.6,3.4);
\node at (2.1,3.9) {$\Delta(s)$};
\node at (-0.95,3.2) {$\xi$};
\draw  (-2.45,3.5) rectangle (-1.45,2.5);
\node at (-1.95,3) {$\frac{\,k_i}{s}$};
\draw[rounded corners]  (-0.45,3.5) rectangle (0.55,2.5);
\node at (0.05,3) {sat$(\cdot)$};
\draw[-latex, line width = .5 pt]  (-1.45,3) -- (-0.45,3);
\draw [-latex, line width = .5 pt] (3.55,3.5) rectangle (4.55,2.5);
\node at (4.05,3) {$\frac{1}{s(s+d)}$};
\draw[-latex, line width = .5 pt] (0.55,3) -- (2.95,3);
\node at (4.8,3.2) {$y$};
\draw [-latex, line width = .5 pt](4.55,3) -- (5,3) -- (5,2) -- (-3.2,2) -- (-3.2,2.95);
\draw [-latex, line width = .5 pt](3.05,3) -- (3.55,3);
\draw [-latex, line width = .5 pt](1.2,3) -- (1.2,3.9) -- (1.6,3.9);
\draw [-latex, line width = .5 pt](2.6,3.9) -- (3,3.9) -- (3,3.05);
\node at (3.2,3.15) {\tiny $+$};
\node at (2.85,2.85) {\tiny $+$};
\draw [-latex](-3.15,3) -- (-2.45,3);
\draw [-latex](-4,3) -- (-3.25,3);
\draw[fill = white] (-3.2,3) circle [radius=0.05]; 
\node at (-3.35,3.15) {\tiny $+$};
\node at (-3.05,2.85) {\tiny $-$};
\node at (-3.65,3.2) {$r$};
\end{tikzpicture}
\centering
\caption{The perturbed system associated with~\eqref{eq:system-MSD}.}
\label{fig:feedback_interconnection_example_CDC_perturbed}
\end{figure}%

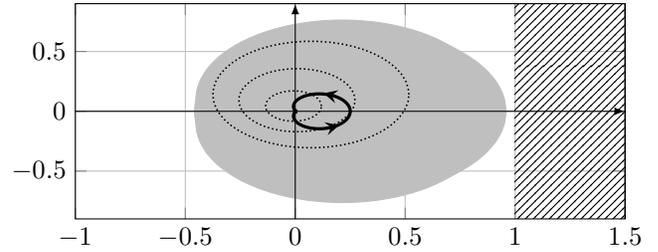
\begin{figure}[t!]
\centering
\begin{tikzpicture}[black]
\begin{axis}[
height = 0.25\textwidth,
width = 0.5\textwidth,
xmin=-1,xmax=1.5,
ymin=-0.9,ymax=0.9,
grid
]
\addplot[pattern=north east lines]  coordinates  {(1, -1.25)  (1.5, -1.25) (1.5, 1.25) (1.0, 1.25) };
\addplot+[name path=A,gray!50,smooth,no markers] table [x index = {0}, y index = {1},  col sep=comma]{msd_cdc_convex_hull.csv};
\addplot+[name path=B, gray!50, no markers]  (0,0); 
\addplot[gray!50] fill between[of=A and B];
\addplot [densely dotted, semithick]  table [x index = {0}, y index = {1},  col sep=comma]{msd_cdc_circle1.csv};
\addplot [densely dotted, semithick]  table [x index = {0}, y index = {1},  col sep=comma]{msd_cdc_circle2.csv};
\addplot [densely dotted, semithick]  table [x index = {0}, y index = {1},  col sep=comma]{msd_cdc_circle3.csv};
\addplot[-latex,black,    domain=-0.9:0.9]  (0,{x}); 
\addplot[-latex, black,    domain=-1:1.5 ]  ({x},0); 
\addplot [black, very thick, %
postaction={decorate}, %
decoration={
  markings,
  mark=at position 0.5*\pgfdecoratedpathlength-10pt
  with {\arrow[xshift=2.5\pgflinewidth,>=stealth]{>}},
  mark=at position 0.5*\pgfdecoratedpathlength+10pt
  with {\arrow[xshift=2.5\pgflinewidth,>=stealth]{>}}
} %
]  table [x index = {0}, y index = {1},  col sep=comma]{msd_cdc_nyquist.csv};
\end{axis}
\end{tikzpicture}
\centering
\caption{The Nyquist diagram of the $\lambda$-shifted transfer function associated with ${G(s)= \frac{k_i}{s^2(s+d)}}$ (solid) lies to the right of the disk $D(-1,0)$ (diagonal lines) for ${d=5}$, ${k_i=-1}$, ${\lambda =1}$.  Robust strict $2$-dominance is guaranteed when the uncertain Nyquist diagram (shaded) given by the envelope of all circles of center $G_\lambda (i \omega)$ and radius $ \norm{\Delta}_{\H_{\infty,0}(\S_{\Lambda})}  |G_\lambda (i \omega)|$ (dotted) lies outside the disk $D(-1,0)$.
} 
\label{fig:msd_CDC19}
\end{figure}

For illustration, assume $\tau$ is sufficiently small (so that $-\tfrac{1}{\tau}$ is to the left of the strip $\SL$). Rewrite the perturbed dynamics as in 
Fig.~\ref{fig:feedback_interconnection_example_CDC_perturbed} by considering a multiplicative uncertainty $\Delta$ such that  ${1+ \Delta(s) = \frac{1}{1+s\tau}}$, \textit{i.e.}
\beq 
\Delta(s) = -\frac{s\tau }{1+s\tau} .
\eeq
By Theorem~\ref{thm:small-gain},  strict $2$-dominance of the nominal system~\eqref{eq:system-MSD} is preserved if the product of the $\L_{2,0}$-gain of the perturbation $\Delta$  and the $\L_{2,2}$-gain of the nominal system~\eqref{eq:system-MSD}  is less than one. By Theorem~\ref{thm:pgain}, the former can be computed as the $\H_{\infty,0}(\S_{\Lambda})$ of the transfer function $\Delta$  in any strip $\SL$ with $\Lambda \subset (0,d)$; the latter can be computed as the supremum over all solutions  of the LMI~\eqref{eq:linear-LMI} for ${\lambda \in \Lambda}$. For example, let ${d=5}$, ${k_i=-1}$, ${\Lambda=(1,2)}$ and ${\tau = 0.1}$. Then
\begin{align*}
 \norm{\Delta}_{\H_{\infty,0}(\S_{\Lambda})} 
 &= \max\{ 
	\norm{\Delta}_{\L_{\infty}(\Li_{1})}, 
	\norm{\Delta}_{\L_{\infty}(\Li_{2})}
	\} \\
&= \max\{ 
	 1.1111, 
	1.0526
	\} = 1.1111 
\end{align*}
and, by Theorem~\ref{thm:pgain}, the $\L_{2,2}$-gain of system~\eqref{eq:system-MSD} is
\begin{align*}
 \gamma_{\Lambda}
 = \sup_{\lambda \in (1,2)} \gamma_{\lambda} 
= \max\{ 
	 0.3528, 
	 0.1414
	\} 
= 0.3528 .
\end{align*}
 We conclude that the perturbed closed-loop system remains  strictly $2$-dominant with any rate ${\lambda \in  (1,2)}$ and, thus, oscillatory, as the perturbation $\Delta$ preserves the unstable equilibrium at the origin.  Note that
 the perturbed closed-loop system can actually tolerate a perturbation $\Delta$ with ${\norm{\Delta}_{\H_{\infty,0}(\S_{\Lambda})} \approx 2.8345} $  and still preserve strict $2$-dominance (since ${ \gamma_{\Lambda} = 0.3528}$), as illustrated in Fig.~\ref{fig:msd_CDC19}.

\section{Conclusion} \label{sec:conclusion}

The paper has shown that the differential $\L_{2,p}$ gain of a linear, time-invariant, $p$-dominant system is the $\H_{\infty,p}$ norm of its transfer function.  Several  parallels have been  drawn between the classical $\H_{\infty}$ norm and the $\H_{\infty,p}$ norm. This suggests that robust stability  and robust $p$-dominance can be studied along the same lines for linear systems.    Future research should focus on the analysis and design of multistable and oscillatory nonlinear uncertain systems that can be decomposed as the feedback interconnection of a linear, \linebreak time-invariant system with bounded gain uncertainties or nonlinearities. A promising research direction is that of robust $p$-dominance analysis using integral quadratic constraints~\cite{megretski1997system}.

\appendix

\subsection{The bilateral Laplace transform}

This section recalls, for completeness, basic definitions and results related to the Laplace transform~\cite{widder1946laplace,oppenheim1996signals,partington1997interpolation}.

\begin{definition}\cite[p.17]{partington1997interpolation} \label{def:Laplace_transform}
Let ${f:\R\to\Co}$ be a measurable function. The \textit{(bilateral) Laplace transform of $f$ at ${s\in \Co}$} is defined as
\beq \label{eq:laplace-transform-definition}
F(s) = \L\{f\}(s) = \int_{-\infty}^{\infty} f(\tau) e^{-s\tau} d\tau
\eeq
for those ${s = \lambda +i\omega}$ such that  ${\int_{-\infty}^{\infty} |f(\tau)| e^{-\lambda \tau} d\tau < \infty  .}$
\end{definition}

\noindent
A complete characterization of the Laplace transform of a function $f$ requires the specification of a \textit{region of convergence}, \textit{i.e.} a set of values ${s\in\Co}$ for which the integral~\eqref{eq:laplace-transform-definition} converges~\cite[p.662]{oppenheim1996signals}.  In general, there may be multiple regions of convergence and  these are always vertical strips in the complex plane, as a consequence of the following result~\cite[p.238]{widder1946laplace}.

\begin{lemma}\label{lemma:laplace1} Let ${f:\R\to\Co}$ be a measurable function  and let ${\Lambda = (\underline{\lambda}, \overbar{\lambda})}$, with ${ -\infty \le \underline{\lambda} < \overbar{\lambda} \le \infty}$.   If the integral
\beq \label{eq:laplace-integral}
\int_{-\infty}^{\infty} f(\tau) e^{-s\tau} d\tau 
\eeq
converges for  ${ \underline{s} = -\underline{\lambda} +i \underline{\omega}}$ and ${\overbar{s} = -\overbar{\lambda} +i \overbar{\omega}}$, then it converges in the strip $\S_{\Lambda}$.
\end{lemma}

\noindent
Thus the region of convergence of a Laplace transform is in general a vertical strip, which may become a half plane, the entire plane or even (parts of) a single vertical line~\cite[p.238]{widder1946laplace}. 
If the Laplace transform converges for ${s\in\S_{\Lambda}}$, with ${\Lambda = (\underline{\lambda}, \overbar{\lambda})}$, and diverges elsewhere, then $-\underline{\lambda}$ and $-\overbar{\lambda}$ are said to be \textit{abscissae of convergence} and the vertical lines ${\Li_{\underline{\lambda}}}$ and ${\Li_{\overbar{\lambda}}}$ are said to be the corresponding \textit{axes of convergence}.
It is clear that if the integral~\eqref{eq:laplace-integral} converges in a strip ${\S_{\Lambda}}$, then it converges uniformly in any closed bounded region inside the strip which does not intersect the boundary of the strip~\cite[p.240]{widder1946laplace}.  Moreover, if the integral~\eqref{eq:laplace-integral} converges along the line ${\Li_\lambda}$ then the region of convergence will be a strip that includes the line ${\Li_\lambda}$~\cite[p.666]{oppenheim1996signals}.

The representation induced by the Laplace transform is unique, as detailed by the following statement~\cite[p.243]{widder1946laplace}.

\begin{lemma} \label{lemma:laplace2}
If ${f:\R\to\Co}$ and ${g:\R\to\Co}$ are measurable functions in any bounded interval and such that
${ \L\{f\} = \L\{g\} }$ in a common region of convergence, then ${f(t)=g(t)}$ for almost every ${t\in\R}$.
\end{lemma}

\noindent
Lemma~\ref{lemma:laplace2} implies that the region of convergence of the Laplace transform of a function ${f\in\L_2(\R)}$ is given by the \textit{intersection} of the regions of convergence of its\footnote{Every ${f\in \L_{2}(\R)}$ admits a unique additive decomposition of the form ${f = f_{+}+f_{-}}$, with ${f_{+}(t)=0}$ for almost all ${t>0}$ and ${f_{-}(t)=0}$ for almost all ${t<0}$. The functions ${f_{+}}$ and ${f_{-}}$ are referred to as the \textit{causal part} of $f$ and \textit{anticausal part} of $f$, respectively.} causal part ${f_{+}}$ and its anticausal part  ${f_{-}}$. If the intersection is non-empty, then ${F(s) = F_{+}(s)  + F_{-}(s) }$. By contrast, when the regions of convergence of ${f_{+}}$ and ${f_{-}}$ do not intersect the Laplace transform of $f$ is not defined, even if the Laplace transforms of ${f_{+}}$ and ${f_{-}}$ are individually well-defined. This is well illustrated by the following example taken from~\cite[p.668]{oppenheim1996signals}.

\begin{example}
Let ${\alpha \in \R}$ and let ${f: t\mapsto e^{-\alpha |t|}}$. The causal part of $f$ is ${f_{+}: t\mapsto e^{-\alpha t} \one(t)}$ and the anticausal part of $f$ is ${f_{-}: t\mapsto e^{\alpha t}\one(-t)}$, where $\one$ denotes the Heaviside unit step function.  The Laplace transform of $f_{
+}$ is ${F_{+}(s) = \frac{1}{s+\alpha}}$, with region of convergence ${\S_{\Lambda_{+}}}$, with ${\Lambda_{+} = (-\infty, \alpha)}$, and the Laplace transform of ${f_{-}}$ is ${F_{-}(s) = -\frac{1}{s-\alpha}}$, with region of convergence ${\S_{\Lambda_{-}}}$, with ${\Lambda_{-} = (-\alpha,\infty)}$. While the Laplace transform of both the causal part and anticausal part of $f$ individually exist, the corresponding regions of convergence do not intersect if ${\alpha \le 0}$, in which case the Laplace transform  of $f$ is not defined. By contrast, if ${\alpha>0}$, then the Laplace transform of $f$ is 
\beq \nn
F(s) = F_{+}(s) + F_{-}(s) =  \frac{1}{s+\alpha} - \frac{1}{s-\alpha} = -\frac{2\alpha}{s^{2}-\alpha^{2}},
\eeq
with region of convergence ${\S_{\Lambda}}$, with ${\Lambda = (-\alpha,\alpha)}$.
\end{example}

The inverse bilateral Laplace transform can be defined using the following result~\cite[p.241]{widder1946laplace}.

\begin{lemma}
Let ${f:\R\to\Co}$ be a measurable function in any bounded interval. Assume that the integral~\eqref{eq:laplace-integral} converges absolutely on the vertical line $\Li_\lambda$ and that $f$ is of bounded variation in a neighbourhood of ${t \in \R}$. Then\footnote{The convention ${f(t) = \frac{1}{2} \lim_{\tau \to t^+} f(\tau) + \frac{1}{2} \lim_{\tau \to t^-} f(\tau)  }$ is used if ${t\in\R}$ is a point of discontinuity of $f$.}
\beq \label{eq:inverse-Laplace-transform}
f(t) = \lim_{\omega \to \infty} \frac{1}{2\pi i}  \int_{-\lambda -i \omega }^{-\lambda +i \omega } F(s) e^{st} ds.
\eeq
\end{lemma}

\noindent
In general, computing inverse Laplace transforms via~\eqref{eq:inverse-Laplace-transform} requires  complex contour integration. In practice, this is often performed using the residue theorem~\cite[p.108]{conway1973functions}.

We conclude this digression on the bilateral Laplace transform with a few words about functions with a \textit{rational} Laplace transform. Definition~\ref{def:Laplace_transform} implies that the region of convergence cannot contain any pole. As a result, if the Laplace transform $F$ of a function $f$ is rational,   its region of convergence is bounded by poles or extends to infinity~\cite[p.669]{oppenheim1996signals}. In particular, the region of convergence is the half plane to the right (left) of the rightmost (leftmost) pole if the anticausal (causal) part of $f$ is zero almost everywhere. 

As a consequence of the residue theorem, the inverse Laplace transform of a rational function can be computed by evaluating~\eqref{eq:inverse-Laplace-transform} via partial fraction expansion and then by inverting each individual term~\cite[p.671]{oppenheim1996signals}. For example, consider the rational function  ${F(s) = \frac{n(s)}{d(s)}},$ where ${n(s) = \sum_{k=0}^{m} n_k s^{k}}$ and ${d(s) = \sum_{k=0}^{n} d_k s^{k}}$, with ${d_{n} = 1}$ and ${n>m}$. If the roots ${\alpha_1, \ldots, \alpha_n}$ of the denominator $d$ are distinct then 
${F(s) = \sum_{k=0}^{n} \frac{f_k}{(s-\alpha_k)},}$
with ${f_k = \frac{n(\alpha_k)}{d^{\prime}(\alpha_k) (s-\alpha_k)}}$. The inverse Laplace transform of each term ${F_k(s) = \frac{f_k}{(s-\alpha_k)}}$ is then determined as follows: if the region of convergence is to the right of the pole  ${s=-\alpha_k}$, then ${\L^{-1}\{F_k\}(t) = f_ke^{-\alpha_k t } \one(t)}$;  if the region of convergence is to the left of the pole  ${s=-\alpha_k}$, then ${\L^{-1}\{F_k\}(t) =-f_ke^{\alpha_k t} \one(-t)}$. A similar argument can be used to find the inverse Laplace transform of a rational function with multiple poles~\cite[p.22]{davies1978integral}.

\end{document}